\documentclass[11pt,a4paper,twoside]{amsart}

\usepackage{enumitem}
\usepackage{amsthm}
\usepackage{bm}
\usepackage{graphics}
\usepackage[mathcal]{euscript}
\usepackage{tikz}
\usepackage{pgfplots}
\pgfplotsset{compat=1.15}
\usetikzlibrary{math, decorations.pathmorphing, patterns}
\usepackage{xpatch}
\usepackage{etoolbox}
\usepackage{xassoccnt}
\usepackage{thmtools}
\usepackage{thm-restate}
\usepackage[colorlinks=true]{hyperref}
\usepackage[capitalise,nameinlink]{cleveref}

\makeatletter
\@namedef{subjclassname@2020}{%
  \textup{2020} Mathematics Subject Classification}
\makeatother

\tikzmath
{
	\simplifyfigures = 0; 
}

\makeatletter
\newcounter{globallistcounter}
\let\realItem\item
\NewDocumentCommand\myItem{o}{%
   \IfNoValueTF{#1}{%
      \realItem%
      \edef\countername{\LastRefSteppedCounter}%
      \def\thegloballistcounter{\csname the\countername\endcsname}%
      \refstepcounter{globallistcounter}%
   }{%
      \realItem[#1]%
      \def\thegloballistcounter{#1}%
      \refstepcounter{globallistcounter}%
   }%
}

\AtBeginEnvironment{enumerate}{%
    \let\item\myItem%
    \stepcounter{globallistcounter}%
}

\def\localItemStringLiteral{LOCALITEM@}
\NewDocumentCommand\myLabel{m}{%
	  \realLabel{\localItemStringLiteral#1}%
      \addtocounter{globallistcounter}{-1}%
      \let\oldthegloballistcounter\thegloballistcounter%
      \def\thegloballistcounter{\theparentcounter\oldthegloballistcounter}%
      \refstepcounter{globallistcounter}%
      \realLabel[\parentcountername]{#1}
}
\NewDocumentCommand\reflocal{m}{%
    \ref{\localItemStringLiteral#1}%
}

\newlist{thmlist}{enumerate}{1}
\setlist[thmlist]{
    before=\let\item\myItem%
           \let\realLabel\label%
           \let\label\myLabel%
           \stepcounter{globallistcounter}%
           \edef\parentcountername{\LastRefSteppedCounter}%
           \edef\theparentcounter{\csname the\parentcountername\endcsname},
    label=\upshape{(\alph{thmlisti})},
    noitemsep}

\declaretheorem[name=Theorem,
                Refname={Theorem,Theorems},
                numberwithin=section]{theo}
\declaretheorem[name=Proposition,
               Refname={Proposition,Propositions},
               numberlike=theo]{prop}
\declaretheorem[name=Lemma,
               Refname={Lemma,Lemmas},
               numberlike=theo]{lemma}

\declaretheorem[name=Definition,
               Refname={Definition,Definitions},
               numberlike=theo]{defi}
\declaretheorem[name=Corollary,
               Refname={Corollary,Corollaries},
               numberlike=theo]{coro}

\declaretheorem[name=Question,
               Refname={Question,Questions},
               numberlike=theo,style=plain]{quest}

\makeatother

\AtEndEnvironment{proof}{\setcounter{claim}{0}}

\newcommand{\bbN}{{\mathbb N}}
\newcommand{\bbR}{{\mathbb R}}

\newcommand{\functional}{F}

\newcommand{\calU}{{\mathcal U}}

\newcommand{\Deltax}{{\Delta x}}
\newcommand{\Deltay}{{\Delta y}}

\newcommand{\eps}{\varepsilon}

\DeclareMathOperator{\sgn}{sgn}

\DeclareMathOperator{\diam}{diam}

\providecommand{\abs}[1]{\lvert#1\rvert}

\newcommand{\smooth}[1]{\tilde{#1}}
\newcommand{\smoothz}[2]{\smooth{#1}_{#2}}
\newcommand{\smoothzp}[3]{\smoothz{#1}{#2,#3}}

\newcommand{\dd}[2]{\frac{\partial #1}{\partial #2}}
\newcommand{\ddxddy}[1]{\dd{#1}{x}\dd{#1}{y}}
\newcommand{\tdd}[2]{\tfrac{\partial #1}{\partial #2}}
\newcommand{\tddxddy}[1]{\tdd{#1}{x}\tdd{#1}{y}}

\begin{document}

\title{Continuity of limit surfaces of locally uniform random permutations}
\author{Jonas Sj{\"o}strand}
\address{School of Education, Culture and Communication, Division of Mathematics and Physics, M\"alardalen University, Box~883, 721~23 V\"aster\aa s, Sweden}
\email{jonas.sjostrand@mdu.se}

\date{September, 2023}

\begin{abstract}
A locally uniform random permutation is generated by sampling $n$ points
independently from some absolutely continuous distribution $\rho$
on the plane and interpreting them as a permutation
by the rule that $i$ maps to $j$ if the $i$th point from the left is the $j$th point from below.
As $n$ tends to infinity, decreasing subsequences in the permutation will appear as curves in the plane,
and by interpreting these as level curves, a union of decreasing subsequences gives rise to a surface. In a recent paper by the author it was shown that, for any $r\ge0$, under the correct scaling as $n$ tends to infinity, the surface of the largest union of $\lfloor r\sqrt{n}\rfloor$ decreasing subsequences approaches a limit in the sense that it will come close to a maximizer of a specific variational integral
(and, under reasonable assumptions, that the maximizer is essentially unique).
In the present paper we show that there exists a continuous maximizer, provided that $\rho$ has bounded density and support.

The key ingredient in the proof is a new theorem about real functions of two variables that are increasing in both variables:
We show that, for any constant $C$,
any such function can be made continuous without increasing the diameter of its image or decreasing anywhere the product of its partial derivatives clipped
by $C$, that is the minimum of the product and $C$.
\end{abstract}

\maketitle

\section{Introduction}
Let $\sigma$ be a \emph{permutation} of $\{1,2,\dotsc,n\}$.
A \emph{subsequence} of $\sigma$ is an ordered sequence
$(\sigma(i_1),\sigma(i_2),\dotsc,\sigma(i_k))$ where
$i_1<i_2<\dotsb<i_k$. It is \emph{increasing} if $\sigma(i_1)<\sigma(i_2)<\dotsb<\sigma(i_k)$
and \emph{decreasing} if $\sigma(i_1)>\sigma(i_2)>\dotsb>\sigma(i_k)$.

The study of increasing and decreasing subsequences in random permutations has
a long history starting with Ulam \cite{Ulam} in 1961, and we refer to Romik \cite{RomikBook} for a comprehensive review.
It has been known since the 1970s by the work of Vershik and Kerov \cite{VershikKerov} that
the longest decreasing (or increasing) subsequence of a random
permutation of $\{1,2,\dotsc,n\}$ has length approximately
$2\sqrt{n}$ for large $n$. More generally, the (scaled) limit
of the cardinality of the largest union of
$\lfloor r\sqrt{n}\rfloor$ disjoint decreasing subsequences is known
for any $r\ge0$, where $\lfloor\cdot\rfloor$ denotes the integral part. This result
follows from Greene's characterization \cite{Greene}
of the shape of the Robinson--Schensted
Young tableaux together with the limit shape found by Vershik and Kerov \cite{VershikKerov}
and independently by Logan and Shepp \cite{LoganShepp}.
A recent paper by the author~\cite{SjostrandIncreasing23} took a step further and looked at the
\emph{location} of these unions. Also, it considered random permutations
not sampled from the uniform distribution but from any instance of a large class of distributions called \emph{locally uniform}.

To be able to define the ``location'' of monotone subsequences, we will consider
permutations embedded in the plane.
Let $\sigma$ be a finite set of points in the plane,
no two of which having the same $x$- or $y$-coordinate.
We can interpret any such $\sigma$ as a permutation
by letting
$\sigma(i)=j$ if the $i$th point from the left is the $j$th point from
below. If $\sigma$ consists of $n$ points that are sampled independently
from some given absolutely continuous distribution $\rho$ on the plane,
$\sigma$ is said to be \emph{locally uniform (with density $\rho$)}.
In particular, if $\rho$ is the uniform distribution on the unit square $(0,1)^2$, then, as a permutation,
$\sigma$ is uniformly distributed among
all permutations of order $n$.

\begin{figure}
\begin{center}
\begin{tikzpicture}
\node[xscale=1,yscale=-1,inner sep=0,outer sep=0](0,0){\includegraphics[scale=0.21]{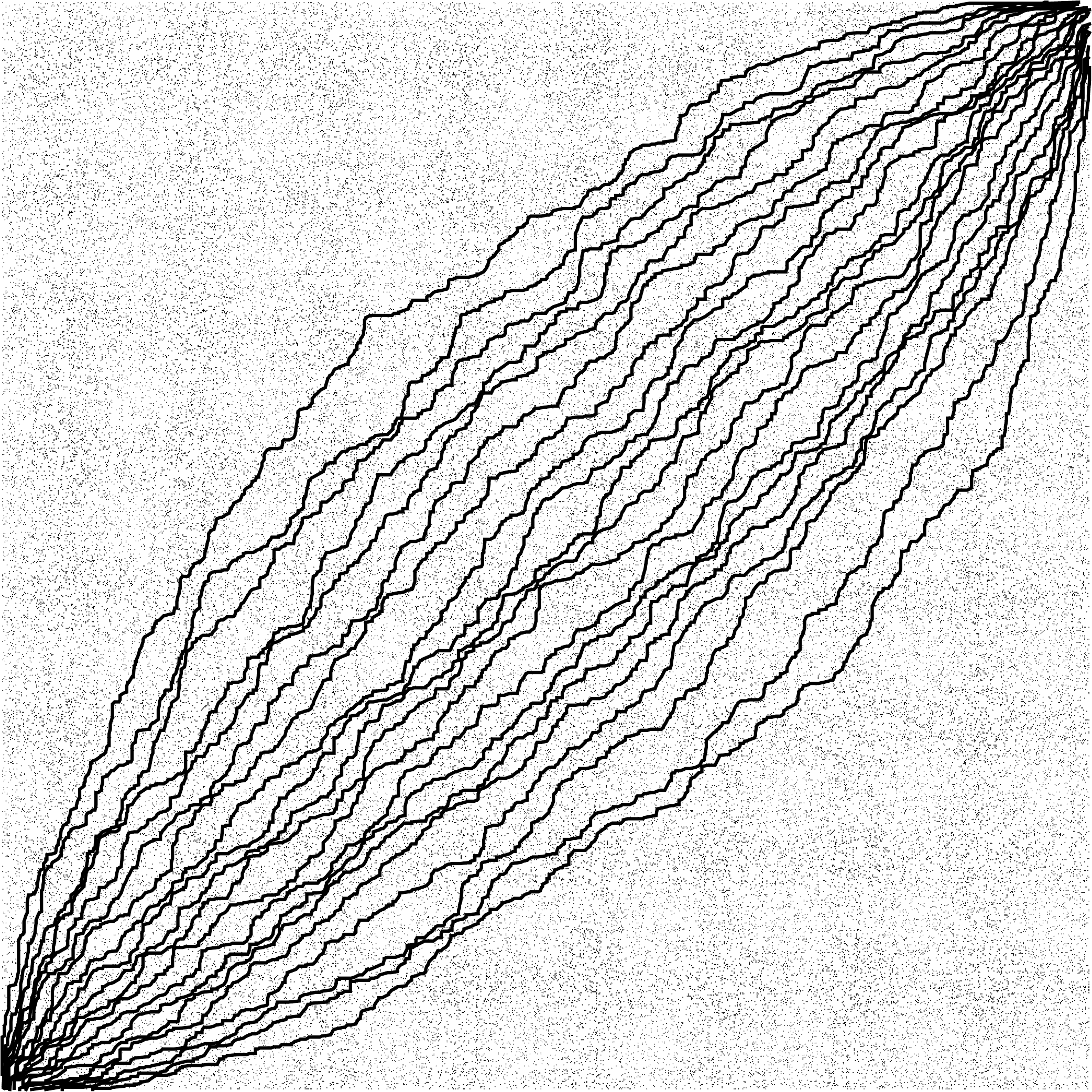}};
\end{tikzpicture}
\end{center}
\caption{The location of the largest union of 20 decreasing subsequences in
a random permutation of order 100{,}000 drawn from the uniform distribution
on the square $(0,1)^2$.
The small dots are the points in the permutation, and adjacent points in
the decreasing subsequences are connected by line segments.}\label{fig:onion}
\end{figure}

In this geometric setting, monotone subsequences of $\sigma$
appear as ``monotone subsets'' of the permutation points in the plane,
and we may speak about the \emph{location} of them.
The location of the longest increasing subsequence was studied by Deuschel and Zeitouni \cite{DeuschelZeitouni95} and by Dauvergne and Vir\'ag \cite{DauvergneVirag21} in the uniform case. 
We will be concerned with the location of a union of decreasing
subsequences, as exemplified in \cref{fig:onion}. Basu et~al.~\cite{BasuEtAl2022}
coined the term \emph{geodesic watermelon} for such a union in the context of
Poissonian last passage percolation.
One could imagine a two-dimensional surface whose level curves follow the
decreasing subsets, and as $n$ tends to infinity, under some rescaling
one might hope to obtain a \emph{limit surface} for a maximal
union of $k$ decreasing subsets, where $k$ depends on $n$.
(It is not hard to see that we must require that $k$ grows as $\sqrt{n}$.)
With the correct definitions, this turns out to be true: In \cite{SjostrandIncreasing23}, it was shown that, for large $n$, the surface
is close to a maximizer of a specific variational integral (and, under reasonable circumstances, that the maximizer is essentially unique).
In the present
paper we study the regularity of such maximizers and show that there is
always a continuous maximizer provided that $\rho$ has bounded density and support.
Note that, in some similar situations, certain surfaces characterized by variational problems turn out to be discontinuous. One example is the recent result of Borga~et~al.~\cite{BorgaEtAl2023determinantal} on
the limit surface of a random standard Young tableux of an L-shape.

Locally uniform random permutations have gained some recent attention:
Dubach \cite{Dubach2023first} studied the length of the
longest increasing subsequence in the case where the density has a singularity. If we replace the assumption that
$\rho$ is absolutely continuous with the weaker assumption that
it has continuous marginals, we obtain something
called \emph{pre-permutons} \cite{Dubach2023first}, and for pre-permutons
that are not absolutely continuous,
Dubach \cite{Dubach2023second} studied the
length of the $k$-th row in the corresponding Young tableaux.
If we forget about the embedding, a pre-permuton $\rho$ induces a probability distribution on $n$-permutations. This distribution is invariant under strictly increasing transformations
of the coordinates, so we can always make the marginal distributions of
$\rho$ uniform on $[0,1]$ without changing the induced distribution on permutations,
that is, we may assume that the support of $\rho$ is contained in $[0,1]^2$
and that $\rho([a,b]\times[0,1])=\rho([0,1]\times[a,b])=b-a$ for any $0\le a\le b\le 1$; see Remark~1.2 in \cite{BorgaEtAl2023}. Such distributions are called
\emph{permutons} \cite{GlebovEtAl2015} and have been studied by several authors,
for instance \cite{HoppenEtAl,Mukherjee2015,BassinoEtAl2022}.

There is a vast literature on the existence and regularity of maximizers or minimizers
$u$ of a variational integral of the form
\[
\int_\Omega f(x,u(x),\nabla u(x))\,dx,
\]
where $\Omega$ is an open subset of $\bbR^m$, $u$ is a function from $\Omega$ to $\bbR^N$,
$\nabla u$ is the total derivative of $u$, and $f$ is a real-valued function on $\Omega\times\bbR^N\times\bbR^{mN}$.
A review of this area of research would lead us too far astray, so
we are content with referring to the two classical books \cite{LadyzhenskayaUraltsevaBook68, MorreyBook68}.
The assumptions may vary regarding
\begin{itemize}
\item
the function space where $u$ lives,
\item
the regularity of $f$ and
\item
the growth conditions for $f$,
\end{itemize}
but typically the admissible $u$ form a Sobolev space, $f$ is supposed to be continuous or even smooth, and there is at least a growth condition of the type
\[
f(x,u,p) \le k-\abs{p}^\alpha
\]
for some constants $\alpha>1$ and $k>0$. None of these assumptions holds in our case. We are interested in monotone functions (in both variables), our $f$ is not required to be continuous in $x$ and it will be nonnegative so it does not satisfy the growth condition above.
There are papers that handle the case where $f$ is discontinuous, like
\cite{CupiniEtAl18, GiaquintaGiusti82}, but to the best of our knowledge no one has
studied the question for monotone functions before.

\section{Terminology and results}\label{sec:results}
A finite set $I$ of points in $\bbR^2$ is \emph{increasing}
if, for any pair of points $(x,y)$ and $(x',y')$ belonging to $I$,
$x<x'$ if and only if $y<y'$. It is \emph{decreasing} if
$x<x'$ if and only if $y>y'$. It is \emph{$k$-increasing}
(resp.~$k$-decreasing) if
it is a union of $k$ increasing (resp.~decreasing) sets.

The following theorem was proved in \cite[Theorem~3.2]{SjostrandIncreasing23}.
\begin{theo}\label{th:narrowrectanglesresult}
Let $\Omega$ be the open rectangle (depicted in \cref{fig:slopedbetaonerectangle})
\[
0<(x+y)/\sqrt2<1,\ \ 0<(y-x)/\sqrt2<\beta
\]
for some $\beta>0$,
and let $r\ge0$.
For each $\gamma>0$, let $\sigma_\gamma$ be a Poisson point
process in the plane with homogeneous intensity $\gamma$.
Define the random variable $\Lambda^{(\gamma)}$ as the size
of a maximal $\lfloor r\sqrt{\gamma}\rfloor$-decreasing subset
of $\sigma_\gamma\cap\Omega$. Then, as $\gamma$
and $\beta$ tends to infinity,
$\Lambda^{(\gamma)}/\beta\gamma$ converges in $L^1$ to a constant $\Phi(r)$.
\end{theo}
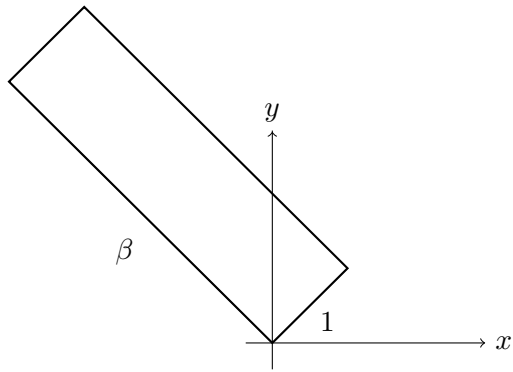
\begin{figure}
\begin{center}
\begin{tikzpicture}[scale=0.7]
\draw [->,thin] (7,-0.5) |- (7,4) node [above] {$y$};
\draw [->,thin] (6.5,0) |- (11,0) node [right] {$x$};
\draw[thick,rotate around={-45:(7,0)}] (0,0) rectangle (7,2);
\draw [rotate around={-45:(7,0)}]
(7,2) -- (7,0) node [black,midway,right,xshift=0pt,yshift=-6pt]
{1};
\draw [rotate around={-45:(7,0)}]
(7,0) -- (0,0) node [black,midway,below,xshift=-6pt,yshift=-6pt]
{$\beta$};
\end{tikzpicture}
\end{center}
\caption{The rectangle $\Omega$ in \cref{th:narrowrectanglesresult}.}
\label{fig:slopedbetaonerectangle}
\end{figure}
The only thing we need to know about the increasing function
$\Phi$ is the following proposition from \cite[Prop.~11.1]{SjostrandIncreasing23}.
\begin{prop}\label{pr:phioneisone}
$\Phi(t)=1$ if $t\ge\sqrt2$.
\end{prop}
Note, however, that is was conjectured in \cite[Conj.~3.5]{SjostrandIncreasing23} that
in fact
\[
\Phi'(r)=
\begin{cases}
\sqrt2-r & \text{if $0\le r \le \sqrt2$,} \\
0 & \text{if $r>\sqrt2$.}
\end{cases}
\]

Define a partial order $\le$
on $\bbR^2$ by letting $(x_1,y_1)\le(x_2,y_2)$ if
$x_1\le x_2$ and $y_1\le y_2$.
For any subset $A$ of $\bbR^2$,
a function $u:A\rightarrow\bbR$ is \emph{doubly increasing} if
$u(x_1,y_1)\le u(x_2,y_2)$ whenever $(x_1,y_1)\le(x_2,y_2)$.
The \emph{diameter} $\sup S-\inf S$
of a bounded subset $S$ of $\bbR$ will be denoted by $\diam S$.
For $r\ge0$, we let $\calU_r(A)$ denote the set of
doubly increasing functions
$u$ on $A$ with $\diam u(A)\le r$, and we let $\calU(A):=\bigcup_{r\ge0}\calU_r(A)$
denote the set of all bounded
doubly increasing functions on $A$.

\begin{figure}
\begin{center}
\begin{tikzpicture}[scale=0.8]
\draw (1,8) -- (1,4) -- (3,4) -- (3,1) -- (8,1);
\draw[fill] (1,4) circle(2pt) (3,1) circle(2pt);
\draw (2,8) -- (2,6) -- (4,6) -- (4,3) -- (5,3) -- (5,2) -- (8,2);
\draw[fill] (2,6) circle(2pt) (4,3) circle(2pt) (5,2) circle(2pt);
\draw (6,8) -- (6,5) -- (8,5);
\draw[fill] (6,5) circle(2pt);
\draw (1,2) node {0} (3,5) node {1} (5,5) node {2} (7,6) node {3};
\end{tikzpicture}
\end{center}
\caption{The $\kappa_P$ function for a set $P$ of six points.}
\label{fig:kappa}
\end{figure}
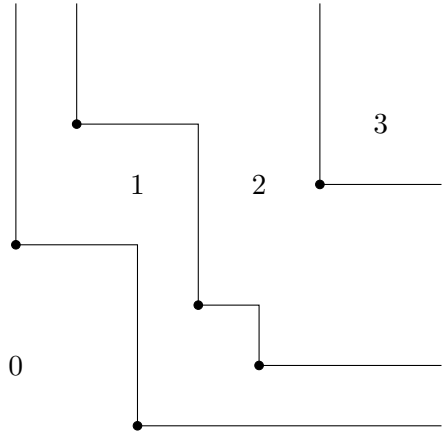

The following well-known combinatorial fact connects decreasing subsets to
increasing ones. A proof was given in \cite{SjostrandIncreasing23}.
\begin{prop}
Let $P$ be a finite set of points in general position in the
sense that no two points share the same $x$- or $y$-coordinate.
Then, $P$ is $k$-decreasing if and only if it has no increasing subset
of cardinality larger than $k$.
\end{prop}
In accordance, our interpretation of point sets as doubly increasing functions looks as follows.
\begin{defi}\label{def:kappa}
For any finite set $P$ of points in the plane,
define a map $\kappa_P:\bbR^2\rightarrow\bbN$ by letting
$\kappa_P(x,y)$ be the maximal size of an increasing subset
of $P\cap((-\infty,x]\times(-\infty,y])$.
\end{defi}
See \cref{fig:kappa} for an example.

Let $\mu$ denote the Lebesgue measure on $\bbR^2$.
By a \emph{bounded probability domain} we will mean a pair $(\Omega,\rho)$ where
$\Omega$ is a bounded open subset of $\bbR^2$
and $\rho$ is a bounded probability density function on $\Omega$,
that is, a bounded nonnegative function on $\Omega$ such that
$\int_\Omega \rho\,d\mu=1$.

\begin{defi}\label{df:FrhoandL}
For any $\eta,\theta\ge0$, let
\[
L(\eta,\theta):=\begin{cases}
\eta\,\Phi(\sqrt{2\theta/\eta}) & \text{if $\eta>0$,} \\
0 & \text{if $\eta=0$,}
\end{cases}
\]
and, for any bounded probability domain $(\Omega,\rho)$, let
$\functional_\rho:\calU(\Omega)\rightarrow\bbR$
be a (nonlinear) functional given by
\[
\functional_\rho(u) := 
\int_\Omega L(\rho, \tddxddy{u})\,d\mu.
\]
\end{defi}

By a result of Burkill and
Haslam-Jones \cite{BurkillHaslamjones}, a doubly increasing
function is differentiable almost everywhere, so the integrand above is defined
almost everywhere.
In \cite[Sec.~4]{SjostrandIncreasing23} it was shown that the integrand is
integrable so that $\functional_\rho$ is well defined.

Now we are in position to state the main result of \cite{SjostrandIncreasing23}. 
\begin{theo}\label{th:limitsurface}
Let $(\Omega,\rho)$ be a bounded probability domain, and,
for each positive integer
$n$, let $\tau_n$ be a set of $n$ random points in $\Omega$ sampled
independently with probability density function $\rho$.
Let $r\ge0$ and $\eps>0$.
Then, with probability tending to one as
$n\rightarrow\infty$, for
every maximal $\lfloor r\sqrt{n}\rfloor$-decreasing subset $P$
of $\tau_n$, the functional $\functional_\rho$ has a
(not necessarily unique) maximizer $u$ over $\calU_r(\Omega)$ such that
$\int_\Omega\abs{\left(\kappa_P/\sqrt{n}\right)-u}\,d\mu<\eps$.
\end{theo}
In fact, as was shown in \cite[Theorem~10.2]{SjostrandIncreasing23}, the above theorem holds even without the assumption of $\Omega$ and $\rho$ being bounded,
but we do not need that level of generality here.

\Cref{th:limitsurface}
guarantees that the functional $\functional_\rho$ has
a maximizer over $\calU_r(\Omega)$, but it does not
say anything about the regularity of the maximizer. Our main result
is the following theorem, stating that there always exists
a continuous maximizer.
\begin{restatable}{theo}{maintheorem}\label{th:main}
Let $(\Omega,\rho)$ be a bounded probability domain.
Then, for any $r\ge 0$, there is a continuous maximizer of $\functional_\rho$ over $\calU_r(\Omega)$.
\end{restatable}
The proof will appear in \cref{sec:mainproof}. Its main ingredient
is the following general theorem about patching the discontinuities of a doubly increasing function, which will be proved in \cref{sec:largeproductproof}.
\begin{restatable}{theo}{largeproducttheorem}\label{th:largeproduct}
Let $u\,:\,\bbR_{\ge0}^2\rightarrow\bbR_{\ge0}$ be a doubly increasing function
such that $u(0,t)=u(t,0)=0$ for all $t\ge0$.
Then, for any $C>0$ there is
a continuous doubly increasing function $\smooth{u}$ on $\bbR_{\ge0}^2$
such that $0\le \smooth{u}\le u$ everywhere and
$\ddxddy{\smooth{u}}\ge C$ on the set where $\smooth{u}\ne u$
and $\smooth{u}$ is differentiable.
\end{restatable}
The rest of the paper is organized as follows.
First, in \cref{sec:righttool} we explain why \cref{th:largeproduct} is the
tool we need for proving our main result, and we show why some of the
seemingly simpler approaches fail.
Then, in \cref{sec:largeproductproof} we prove \cref{th:largeproduct} after warming up
by showing its one-dimensional analogue. 
In \cref{sec:mainproof} we prove our main result, \cref{th:main},
and in \cref{sec:nondifferentiability} we show that this is the best result possible
in the sense that continuity cannot be replaced by differentiability. 
Finally, in \cref{sec:future} we discuss some open problems.

\Cref{sec:mainproof,sec:nondifferentiability,sec:future} are independent of \cref{sec:largeproductproof}, so a reader that
are willing to use \cref{th:largeproduct} without seeing its proof might want to
skip \cref{sec:largeproductproof}.

\section{Is \texorpdfstring{\cref{th:largeproduct}}{this} the right tool for us?}\label{sec:righttool}
To prove \cref{th:main}, we want to show that a doubly increasing function $u$
can be made continuous without increasing the diameter of its image or decreasing the value of $L(\rho,\ddxddy{u})$ anywhere. It is shown in \cite[Lemma~4.9]{SjostrandIncreasing23} that $L$ is increasing in the second variable, so it would be enough if we did not decrease $\ddxddy{u}$ anywhere. Can this be achieved?

To gain some intuition, let us consider the one-dimensional
version of the same problem:
Can an increasing function of a single variable be
made continuous without increasing the diameter of its image
or decreasing its derivative?
The answer is yes. If $u:\bbR_{\ge0}\rightarrow\bbR$ is an increasing function,
let $u'$ be its derivative, which is defined almost everywhere, and define
the increasing function $\smooth{u}:\bbR_{\ge0}\rightarrow\bbR$ by letting
$\smooth{u}(x)=u(0) + \int_0^x u'(t)\,dt$ (in the sense of Lebesgue).
Then $u(0)\le\smooth{u}\le u$ everywhere and $\smooth{u}'=u'$ whenever $u'$ is defined. In terms of the graph of the function, this construction removes a discontinuity by translating the part of the graph that is to the right
of the discontinuity downwards.

For doubly increasing functions, however, the situation is not that
simple. In general its discontinuities cannot be removed by global translations since they might consist of a local slit in the surface or even an isolated point.
So we cannot make $u$ continuous without changing its partial derivatives,
but maybe we can still avoid decreasing $\ddxddy{u}$?

\begin{figure}
\begin{center}
\tikzmath{
	\a = 0.0;
	\d = 0.9;
	\numsamples = ifthenelse(\simplifyfigures,5,21);
} 
\begin{tikzpicture}[scale=3]
	\draw[->] (-1.1, 0) -- (1.1, 0) node[right] {$x$};
	\draw[->] (0, -1.1) -- (0, 1.1) node[above] {$y$};
	\foreach \c in {-1.0, -0.9, ..., 1.0}
	{
		\draw ({(-1-\c)/2}, {(1-\c)/2}) -- ({(1+\c)/2}, {(-1+\c)/2});
    }
    \draw (-0.5,-0.5) node {$-1$};
    \draw (0.5,0.5) node {$+1$};
\end{tikzpicture}
\end{center}
\caption{A level plot of the function $u$ defined by \cref{eq:discontinuityexample}.}
\label{fig:discontinuityexample}
\end{figure}
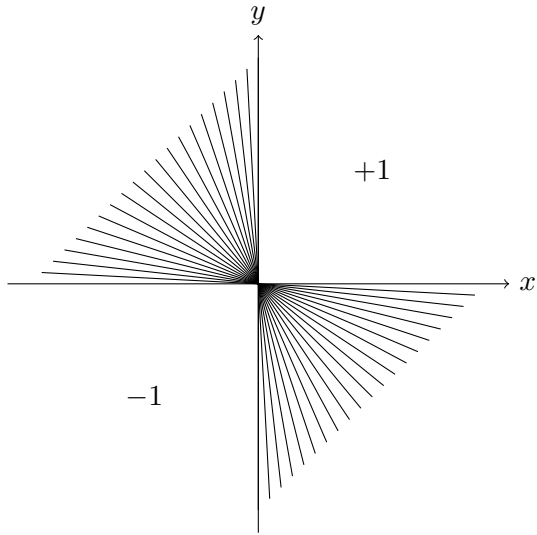

Unfortunately, the answer is no, as the following example shows.
Define the doubly increasing function $u$ on $\bbR^2$ by letting
\begin{equation}\label{eq:discontinuityexample}
u(x,y):=
\begin{cases}
\frac{x+y}{\abs{x-y}} & \text{if $xy\le0$ and $(x,y)\ne(0,0)$,} \\
\sgn(x+y) & \text{if $xy>0$,} \\
0 & \text{if $(x,y)=(0,0).$}
\end{cases}
\end{equation}
As can be seen clearly from its level plot in \cref{fig:discontinuityexample},
$u$ has a discontinuity at the origin and we will show that it is impossible
to get rid of this discontinuity without decreasing the product of the partial
derivatives outside a null set. Suppose $\tilde{u}$ is another doubly increasing
function such that
$\ddxddy{\tilde{u}}\ge\ddxddy{u}$ almost everywhere. For any $\eps>0$,
let $C_\eps$ be the line segment from $(0,-\eps)$ to $(\eps,0)$. A doubly
increasing function is differentiable almost everywhere, so for almost every $\eps>0$
both $u$ and $\tilde{u}$ are differentiable almost everywhere on $C_\eps$.
Consider such an $\eps$, and parameterize $C_\eps$ by
letting $\mathbf{r}(t)=(\eps t,\eps(t-1))$ for $0\le t\le 1$.
Then,
\begin{multline*}
\tilde{u}(\eps,0)-\tilde{u}(0,-\eps)
\ge \int_{C_\eps} \nabla\tilde{u}\cdot d\mathbf{r}
=\int_0^1\left(\tdd{\tilde{u}}{x}(\mathbf{r}(t))+\tdd{\tilde{u}}{y}(\mathbf{r}(t))\right)\eps\,dt \\
\ge2\eps\int_0^1\sqrt{\tdd{\tilde{u}}{x}(\mathbf{r}(t))\tdd{\tilde{u}}{y}(\mathbf{r}(t))}\,dt
\ge2\eps\int_0^1\sqrt{\tdd{u}{x}(\mathbf{r}(t))\tdd{u}{y}(\mathbf{r}(t))}\,dt
\\=4\int_0^1\sqrt{t(1-t)}\,dt
=\frac{\pi}{2},
\end{multline*}
where the second inequality follows from the inequality of arithmetic and
geometric means.
This shows that $\tilde{u}(\eps,0)-\tilde{u}(0,-\eps)$ does not tend to zero as $\eps$ tends
to zero, and we conclude that $\tilde{u}$ is discontinuous at the origin, just like $u$.

We have seen that in general we cannot make a doubly increasing function $u$ continuous without decreasing $\ddxddy{u}$, but by \cref{pr:phioneisone}, as long as $\ddxddy{u}$
is larger than $\rho$ its exact value does not matter since $L(\rho,\ddxddy{u})$ will
be equal to $\rho$ in any case. Thus, in order to prove \cref{th:main}
it seems to be enough to have a result along the lines of \cref{th:largeproduct}, and as we have seen, this result is essentially the best we can hope for.

\section{Making a doubly increasing function continuous}\label{sec:largeproductproof}
The goal of this section is to prove \cref{th:largeproduct}, which will be the main
ingredient in the proof of our main result in \cref{sec:mainproof}.
As a warm-up, we will consider the one-dimensional analogue of \cref{th:largeproduct}.

\subsection{The one-dimensional case}
\begin{theo}\label{th:onedimlargeproduct}
Let $u$ be an increasing function on $\bbR_{\ge0}$ with $u(0)=0$.
Then, for any $a>0$ there is
a continuous increasing function $\smooth{u}$ on $\bbR_{\ge0}$
such that $0\le \smooth{u}\le u$ everywhere and
$\smooth{u}'=a$ on the set where $\smooth{u}\ne u$
and $\smooth{u}$ is differentiable.
\end{theo}
\begin{proof}
Define the function $\smooth{u}\,:\,\bbR_{\ge0}\rightarrow\bbR$ by
\[
\smooth{u}(x) := \inf_{0\le x'\le x} \bigl(u(x')+a(x-x')\bigr).
\]
In words, $\smooth{u}(x)$ is the maximal height $h$ such that
the line with slope $a$ passing through the point $(x,h)$ lies
weakly below the graph of $u$ on the interval $[0,x]$. Our goal is
to show that $\smooth{u}$ has the properties stated in the theorem.

It follows directly from the definition that $0\le\smooth{u}\le u$ everywhere.
To show that $\smooth{u}$ is increasing, we reparameterize and write
\[
\smooth{u}(x)=\inf_{\Deltax \ge 0}
\bigl(u(\max\{x-\Deltax,0\})+a\Deltax\bigr).
\]
Since $u$ is increasing, it follows that $\smooth{u}$ is increasing.

For any
$0\le x_1\le x_2$, we have
\begin{multline*}
\smooth{u}(x_2)=\inf_{0\le x'\le x_2} \bigl(u(x')+a(x_2-x')\bigr)
\le \inf_{0\le x'\le x_1} \bigl(u(x')+a(x_2-x')\bigr)\\
=\smooth{u}(x_1) + a(x_2-x_1),
\end{multline*}
and hence $\smooth{u}$ is continuous.

Finally, let $x$ be a point where $\smooth{u}(x)<u(x)$ and where $\smooth{u}$
is differentiable.
Let $\delta:=\bigl(u(x)-\smooth{u}(x)\bigr)/a$.
Then, for any $x'$ in the interval $x\le x'\le x+\delta$, it holds
that $x-x'\ge-\delta$ and hence $u(x') + a(x-x') \ge u(x) - a\delta = \smooth{u}(x)$. It follows that,
for any $0 \le \Deltax \le \delta$,
\[
\smooth{u}(x)=\inf_{0\le x'\le x+\Deltax}\bigl(u(x') + a(x-x')\bigr).
\]
Also, by definition,
\[
\smooth{u}(x+\Deltax)=\inf_{0\le x'\le x+\Deltax}\bigl(u(x') + a(x+\Deltax-x')\bigr),
\]
so $\smooth{u}(x+\Deltax)-\smooth{u}(x) = a\Deltax$ and it follows that
$\smooth{u}'(x)=a$.
\end{proof}

\subsection{The two-dimensional case}
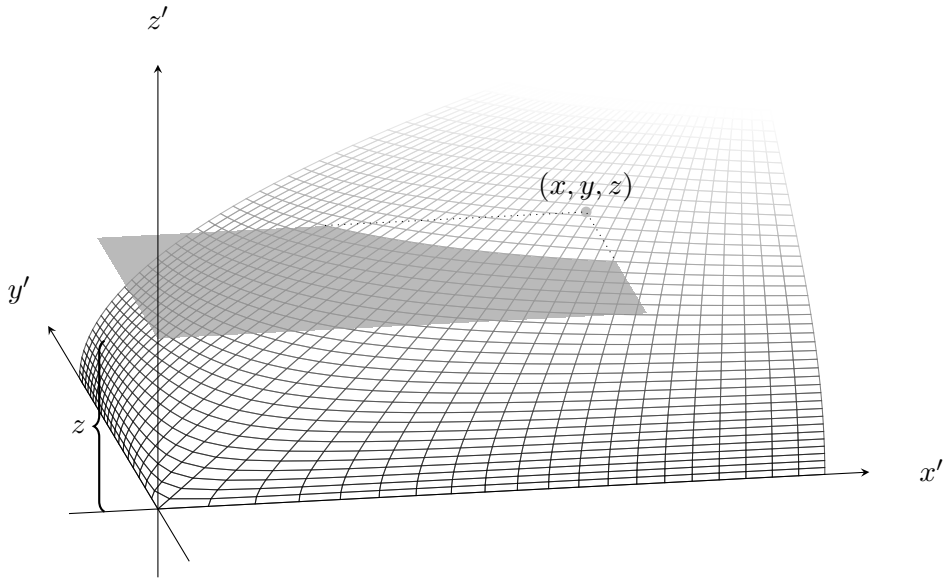
\begin{figure}
\begin{center}
\tikzmath{
	\numsamplesmesh = ifthenelse(\simplifyfigures,9,41);
	\numsamplessurf = ifthenelse(\simplifyfigures,11,21);
} 
\begin{tikzpicture}
[declare function={%
		h(\t,\y) = ((1-\t)/2 + sqrt((1-\t)^2/4 + \y^2*\t);
		f(\x,\y) = (\x > 0) * h(cot(deg(\x)),\y) + (\x <= 0) * \y^2;
		g(\x,\y) = (\x < pi/2-0.0001) * h(tan(deg(\x)),\y) + (\x >= pi/2-0.0001) * \y^2;
		F(\s) = \s + (\s*(1-\s))^2;
		a(\s) = pi/4 + \s * 0.4;	
		qx(\s) = f(a(\s), 0.5);
		qy(\s) = F(g(a(\s), 0.5));
		px(\s) = (\s > 0) * \s * qx(1);
		py(\s) = (\s < 0) * (-\s) * qy(-1);
	}]%
    \begin{axis}[width=\textwidth,
    axis line style={draw=none},
    xmin=-0.2,xmax=1.6,ymin=-0.4,ymax=1.4,zmin=-0.2,zmax=1.3,
    xtick=\empty,ytick=\empty,ztick=\empty,
    axis lines=center,
    colormap/blackwhite,
    view={-10}{25}]

    \addplot3[mesh,shader=interp,opacity=0.7,domain=0:pi/2,y domain=0:0.8,
    samples=\numsamplesmesh]
		({1.5*f(x,y)}, {F(g(x,y))}, y);

    \addplot3[surf,shader=interp,opacity=0.7,domain=-1:1,y domain=0:1,%
    samples=\numsamplessurf]
		({1.5*((1-y)*px(x) + y*qx(x))}, {(1-y)*py(x) + y*qy(x)}, 0.5);

	\fill[gray!50] ({1.5*qx(1)},{qy(-1)},0.5) circle (2pt);
	\draw ({1.5*qx(1)},{qy(-1)},0.5) node[anchor=south]{$(x,y,z)$};
	\draw[dotted] ({1.5*qx(-1)},{qy(-1)},0.5) -- ({1.5*qx(1)},{qy(-1)}, 0.5) -- ({1.5*qx(1)},{qy(1)},0.5);
	\draw [thick,decorate,decoration={brace,amplitude=4pt},xshift=0cm,yshift=0pt]
      (-0.12,0,0) -- (-0.12,0,0.5) node [midway,left,xshift=-0.1cm] {$z$};

    \end{axis}

    \begin{axis}[width=\textwidth,
    xmin=-0.2,xmax=1.6,ymin=-0.4,ymax=1.4,zmin=-0.2,zmax=1.3,
    xtick=\empty,ytick=\empty,ztick=\empty,
    axis lines=center,
    colormap/blackwhite,
    view={-10}{25},
    xlabel={$x'$},
    ylabel={$y'$},
    zlabel={$z'$},
    xlabel style={at={(ticklabel* cs:1.05)},anchor=west},
    ylabel style={at={(ticklabel* cs:1.05)},anchor=south east},
    zlabel style={at={(ticklabel* cs:1.05)},anchor=south}]
    \end{axis}
\end{tikzpicture}
\end{center}
\caption{An $(x,y,z)$-ceiling (shaded) of a surface $z'=u(x',y')$ (meshed).}
\label{fig:xyzceiling}
\end{figure}

Let $\calU$ be the set of doubly increasing functions $u$ on
$\bbR_{\ge0}^2$ such that $u(0,t)=u(t,0)=0$ for any $t\ge0$.
Fix any positive $a$. (Later on, in the proof of \cref{th:largeproduct}, we
will put $a=4C$.) For any $(x,y,z)\in\bbR_{\ge0}^3$, the horizontal
planar region $\{(x',y',z)\ :\ 0\le x'\le x,\ 0\le y'\le y,\ u(x',y')\le z\}$
is called the $(x,y,z)$-ceiling of $u$; see \cref{fig:xyzceiling} for an illustration.
A plane
in $\bbR^3$ is called an $a$-plane if the product of its
x-slope and y-slope is $a$.
\begin{defi}\label{df:smoothuzp}
For any $u$ in $\calU$ and any $z\ge0$ and $p>0$, let
\[
\smoothzp{u}{z}{p}(x,y):=
\inf_{
\substack{
0\le x'\le x\\
0\le y'\le y\\
u(x',y')\le z}}
\bigl(z+p(x-x')+\tfrac{a}{p}(y-y')\bigr).
\]
\end{defi}
The $a$-plane with x-slope $p$ passing through the point $(x,y,h)$ lies
weakly below a point $(x',y',z)$ if and only if $h\le z+p(x-x')+\tfrac{a}{p}(y-y')$.
Thus, $\smoothzp{u}{z}{p}(x,y)$ is the maximal height $h$ such that
the $a$-plane with x-slope $p$ passing through the point $(x,y,h)$ lies
weakly below the $(x,y,z)$-ceiling of $u$.
\begin{defi}\label{df:smoothuz}
For any $u$ in $\calU$ and any $z\ge0$, let
\[
\smoothz{u}{z}(x,y):=\sup_{p>0}\smoothzp{u}{z}{p}(x,y).
\]
\end{defi}
We will show in \cref{lm:pattains} that the supremum is attained for some $p$,
so we can think of $\smoothz{u}{z}(x,y)$ as the maximal height $h$ such that
there is an $a$-plane through the point $(x,y,h)$ that lies
weakly below the $(x,y,z)$-ceiling of $u$.

Imagine that $u$ had the property that for each point $(x,y,u(x,y))$ and each
$z<u(x,y)$ there is some $a$-plane through $(x,y,u(x,y))$ that lies
weakly below the $(x,y,z)$-ceiling. Intuitively, that would essentially mean that
from each point on the $u$-surface, all points to the south-east
are within reach without requiring us to climb down too steeply. If $u$ does not have this property,
the idea is to lower the faulting points just enough to enforce the property.
In other words, we want the surface to be as high as possible at $(x,y)$,
but low enough to guarantee that
any lower $(x,y,z)$-ceiling can be reached along some $a$-plane.
This motivates our final definition.
\begin{defi}\label{df:smoothu}
For any $u$ in $\calU$, let
\[
\smooth{u}(x,y):=\inf_{z\ge0}\smoothz{u}{z}(x,y).
\] 
\end{defi}

Our goal is to show that, for large enough $a$,
$\smooth{u}$ has the properties stated in \cref{th:largeproduct}.
First we present two technical lemmas.

\begin{lemma}\label{lm:infdiffbounds}
For any real-valued functions $f$ and $g$ defined on the same set $X$, it holds that
\[
\inf_{x\in X}\bigl(f(x)-g(x)\bigr) \le \inf_{x\in X} f(x) - \inf_{x\in X} g(x)
\le\sup_{x\in X}\bigl(f(x)-g(x)\bigr).
\]
\end{lemma}
\begin{proof}
Since infimum of functions is a superadditive operation, we have
\[
\inf f=\inf(f-g+g) \ge \inf(f-g) + \inf g,
\]
which proves the first inequality of the lemma. By symmetry, we may exchange $f$ and $g$ in the
proven inequality and obtain
\[
\inf g - \inf f \ge \inf(g-f) = -\sup(f-g),
\]
which proves the second inequality of the lemma.
\end{proof}

\begin{lemma}\label{lm:infequicontinuous}
Let $\mathcal{F}$ be a family of real functions defined on some interval $I\subseteq\bbR$.
Suppose $\mathcal{F}$ is pointwise equicontinuous and bounded below in the sense that
$g(x):=\inf_{f\in\mathcal{F}} f(x)$ exists for any $x\in I$. Then the function
$g:I\rightarrow\bbR$ is continuous.
\end{lemma}
\begin{proof}
Take any $x\in I$ and any $\delta>0$. By equicontinuity, there is an
$\eps>0$ such that $\abs{f(x)-f(y)}<\delta$ for any $f\in\mathcal{F}$ and any
$y\in I$ with $\abs{x-y}<\eps$.
Take such an $\eps$ and take any $y\in I$ with $\abs{x-y}<\eps$. 
From \cref{lm:infdiffbounds} applied to the functions $f\mapsto f(x)$ and
$f\mapsto f(y)$ we obtain
\[
-\delta
\le
\inf_{f\in\mathcal{F}}\bigl(f(x)-f(y)\bigr)
\le
\inf_{f\in\mathcal{F}}f(x) - \inf_{f\in\mathcal{F}}f(y)
\le
\sup_{f\in\mathcal{F}}\bigl(f(x)-f(y)\bigr)
\le
\delta.
\]
It follows that $\abs{g(x)-g(y)}\le\delta$ and thus that $g$ is continuous.
\end{proof}

The next lemma collects some basic properties of the functions $\smoothzp{u}{z}{p}$,
$\smoothz{u}{z}$ and $\smooth{u}$.
\begin{lemma}\label{lm:basicsmoothing}
The following holds.
\begin{thmlist}
\item\label{lm:basicsmoothingA}
$\smoothz{u}{z}(x,y)\ge\smoothzp{u}{z}{p}(x,y)\ge z$ with equalities if $z\ge u(x,y)$.
\item\label{lm:basicsmoothingB}
$0\le\smooth{u}\le u$.
\item\label{lm:basicsmoothingC}
$\smoothzp{u}{z}{p}$, $\smoothz{u}{z}$ and $\smooth{u}$ are all
doubly increasing.
\item\label{lm:basicsmoothingD}
$\smooth{u}(t,0)=\smooth{u}(0,t)=0$ for any $t\ge0$.
\end{thmlist}
\end{lemma}
\begin{proof}
\reflocal{lm:basicsmoothingA} follows directly from the definitions.
To prove \reflocal{lm:basicsmoothingB} we put $z:=u(x,y)$ and
obtain from \reflocal{lm:basicsmoothingA} that
$\smoothz{u}{z}(x,y)=z$ and hence $\smooth{u}(x,y)\le z$.
For \reflocal{lm:basicsmoothingC} it suffices to show that
$\smoothzp{u}{z}{p}$ is doubly increasing. Reparameterizing \cref{df:smoothuzp}
we can write
\[
\smoothzp{u}{z}{p}(x,y)=\inf_{(\Deltax,\Deltay) \in S_{x,y}}
\bigl(z+p\Deltax+\tfrac{a}{p}\Deltay\bigr),
\]
where $S_{x,y}:=\{(\Deltax,\Deltay)\in\bbR_{\ge0}^2\,:\,u(\max\{x-\Deltax,0\},\max\{y-\Deltay,0\})\le z\}$.
Since $u$ is doubly increasing, for any $(x_1,y_1)\le(x_2,y_2)$
we have $S_{x_1,y_1}\supseteq S_{x_2,y_2}$ and it follows that 
$\smoothzp{u}{z}{p}$ is doubly increasing.

Finally, \reflocal{lm:basicsmoothingD} follows directly from \reflocal{lm:basicsmoothingB}.
\end{proof}

To show that the supremum in \cref{df:smoothuz} is attained, we need to show that
$\smoothzp{u}{z}{p}$ varies continuously with $p$.
\begin{lemma}\label{lm:pcontinuous}
For any fixed $x,y,z\ge0$, $\smoothzp{u}{z}{p}(x,y)$ is continuous
as a function of $p$.
\end{lemma}
\begin{proof}
The functions $f_{x',y'}(p):=z+p(x-x')+\frac{a}{p}(y-y')$
for $0\le x'\le x$, $0\le y'\le y$ are pointwise equicontinuous since
$x-x'$ and $y-y'$ are bounded by $x$ and $y$, respectively.
Hence, by \cref{lm:infequicontinuous},
\[
\smoothzp{u}{z}{p}(x,y)=
\inf_{
\substack{
0\le x'\le x\\
0\le y'\le y\\
u(x',y')\le z}}
f_{x',y'}(p)
\]
is continuous in $p$.
\end{proof}

\begin{lemma}\label{lm:pattains}
For each $x,y,z\ge0$ there is a $p>0$ such that
$\smoothz{u}{z}(x,y)=\smoothzp{u}{z}{p}(x,y)$.
\end{lemma}
\begin{proof}
Let $(p_n)$ be a sequence of positive numbers.
By \cref{lm:basicsmoothingA}, we have $\smoothzp{u}{z}{p_n}(x,y)\ge z$.
If $p_n\rightarrow 0$ it holds that $\smoothzp{u}{z}{p_n}(x,y)\rightarrow z$
since $z+p_n(x-0)+\frac{a}{p_n}(y-y)\rightarrow z$.
If $p_n\rightarrow \infty$ we have $\smoothzp{u}{z}{p_n}(x,y)\rightarrow z$
since $z+p_n(x-x)+\frac{a}{p_n}(y-0)\rightarrow z$.
Since $\smoothzp{u}{z}{p}(x,y)\ge z$ and $\smoothzp{u}{z}{p}(x,y)$ is
continuous as a function of $p$ (by \cref{lm:pcontinuous}),
that function attains its supremum.
\end{proof}

Now we are ready to prove that $\smooth{u}$ is continuous.
\begin{lemma}\label{lm:smoothcontinuous}
For any $0\le x'\le x$, $0\le y'\le y$, we have
\begin{equation}\label{eq:usmalldifference}
\smooth{u}(x,y)-\smooth{u}(x',y')\le
\sqrt{a\bigl(y(x-x')+x(y-y')\bigr)}.
\end{equation}
In particular, $\smooth{u}$ is continuous.
\end{lemma}
\begin{proof}
Take any $0\le x'\le x$, $0\le y'\le y$ and let $k:=\sqrt{a\bigl(y(x-x')+x(y-y')\bigr)}$. If $k=0$, either $(x,y)=(x',y')$ and \cref{eq:usmalldifference} follows trivially, or
$x=x'=0$ or $y=y'=0$ in which case \cref{eq:usmalldifference} follows
from \cref{lm:basicsmoothingD}. Thus, we may assume that $k>0$.

By the second inequality of \cref{lm:infdiffbounds} applied to the functions $z\mapsto \smoothz{u}{z}(x,y)$ and
$z\mapsto \smoothz{u}{z}(x',y')$,
it suffices to show that
\begin{equation}\label{eq:smalldifference}
\smoothz{u}{z}(x,y)-\smoothz{u}{z}(x',y')\le k
\end{equation}
for any $z\ge 0$.
If $\smoothz{u}{z}(x,y)\le z+k$, \cref{eq:smalldifference} follows from the fact that
$\smoothz{u}{z}(x',y')\ge z$ (\cref{lm:basicsmoothingA}), so in the following we assume that
$\smoothz{u}{z}(x,y)> z+k$. By \cref{lm:pattains},
there is a $p>0$ such that $\smoothz{u}{z}(x,y)=\smoothzp{u}{z}{p}(x,y)$.
That means that $z+p(x-x'')+\frac{a}{p}(y-y'')>z+k$
for any $0\le x''\le x$ and $0\le y''\le y$ such that
$u(x'',y'')\le z$. Choosing $(x'',y'')=(0,y)$ we obtain
$px>k$, and instead choosing $(x'',y'')=(x,0)$ we
obtain $\frac{a}{p}y>k$. Inverting these inequalities
yields
\[
\frac{a}{p}<\frac{ax}{k}\text{\ \ \ and\ \ \ }
p<\frac{ay}{k}.
\]
Thus, for any $0\le x''\le x'$ and $0\le y''\le y'$, we have
\begin{multline*}
\bigl[z+p(x-x'')+\tfrac{a}{p}(y-y'')\bigr]
-\bigl[z+p(x'-x'')+\tfrac{a}{p}(y'-y'')\bigr]=
\\=p(x-x')+\tfrac{a}{p}(y-y')<[ay(x-x')+ax(y-y')]/k=
k.
\end{multline*}
By the second inequality of \cref{lm:infdiffbounds} applied to the functions
$(x'',y'')\mapsto z+p(x-x'')+\frac{a}{p}(y-y'')$ and $(x'',y'')\mapsto z+p(x'-x'')+\frac{a}{p}(y'-y'')$, it follows that
$\smoothzp{u}{z}{p}(x,y)-\smoothzp{u}{z}{p}(x',y')\le k$,
which implies \cref{eq:smalldifference}.
\end{proof}

Our next lemma bounds the local variation of $\smoothz{u}{z}$ from below.
\begin{lemma}\label{lm:largedifference}
Let $x,y,z\ge 0$ and suppose $u(x,y)>z$. Then
$\smoothz{u}{z}(x+\Deltax,y+\Deltay)-\smoothz{u}{z}(x,y)
\ge\sqrt{a\Deltax\Deltay}$ for any $\Deltax,\Deltay>0$.
\end{lemma}
\begin{proof}
By \cref{lm:pattains}, there is a $p>0$ such that
$\smoothz{u}{z}(x,y)=\smoothzp{u}{z}{p}(x,y)$. Let
$w:=\smoothz{u}{z}(x,y)-z$ and
\[
q:=p\sqrt{\frac{w+\frac{a}{p}\Deltay}{w+p\Deltax}}.
\]
This is well defined since $w\ge0$ by \cref{lm:basicsmoothingA}.
For notational convenience, define the functions
\begin{align*}
f(x',y')&:=z+p(x-x')+\tfrac{a}{p}(y-y')\text{\ and}\\
g(x',y')&:=z+q(x+\Deltax-x')+\tfrac{a}{q}(y+\Deltay-y'),
\end{align*}
and let
\[
A:=\{(x',y')\ :\ 0\le x'\le x+\Deltax,\ 0\le y'\le y+\Deltay,\
u(x',y')\le z\}
\]
be the $(x+\Deltax,y+\Deltay,z)$-ceiling of $u$.
Our goal is to show that
$g(x',y')\ge\smoothz{u}{z}(x,y)+\sqrt{a\Deltax\Deltay}$ for
any $(x',y')\in A$.

Let $S_1$ be the closed line segment from $(x,y-\frac{wp}{a})$
to $(x-\frac{w}{p},y)$ and let $S_2$ be the closed line segment from
$(x+\Deltax,y-\frac{wp}{a})$ to $(x-\frac{w}{p},y+\Deltay)$;
see \cref{fig:triangles}.
\begin{figure}
\begin{center}
\begin{tikzpicture}
\tikzmath{
	\dx = 4;
	\dy = 1;
	\woverp = 3; 
	\wpovera = 2.3;
	\x = 0;
	\y = 0;
} 
\filldraw (\x,\y) circle (2pt) node[anchor=south]{$(x,y)$};
\filldraw (\x-\woverp,\y) circle (2pt) node[anchor=east]{$(x-\frac{w}{p},y)$};
\filldraw (\x,\y-\wpovera) circle (2pt) node[anchor=north]{$(x,y-\frac{wp}{a})$};
\filldraw (\x-\woverp,\y+\dy) circle (2pt) node[anchor=south east]{$(x-\frac{w}{p},y+\Deltay)$};
\filldraw (\x+\dx,\y-\wpovera) circle (2pt) node[anchor=north]{$(x+\Deltax,y-\frac{wp}{a})$};
\filldraw (\x+\dx,\y+\dy) circle (2pt) node[anchor=south]{$(x+\Deltax,y+\Deltay)$};

\draw (\x-\woverp,\y) -- node[below] {$S_1$} (\x,\y-\wpovera);
\draw (\x-\woverp,\y+\dy) -- node[below] {$S_2$} (\x+\dx,\y-\wpovera);

\draw[dotted] (\x-\woverp,\y) -- (\x,\y) -- (\x,\y-\wpovera);
\draw[dotted] (\x-\woverp,\y+\dy) -- (\x+\dx,\y+\dy) -- (\x+\dx,\y-\wpovera);

\draw (\x-\woverp/3,\y-\wpovera/3) node {$T_1$};
\draw (\x-\woverp/3+2*\dx/3,\y-\wpovera/3+2*\dy/3) node {$T_2$};
\end{tikzpicture}
\end{center}
\caption{The situation in the proof of \cref{lm:largedifference}.}\label{fig:triangles}
\end{figure}
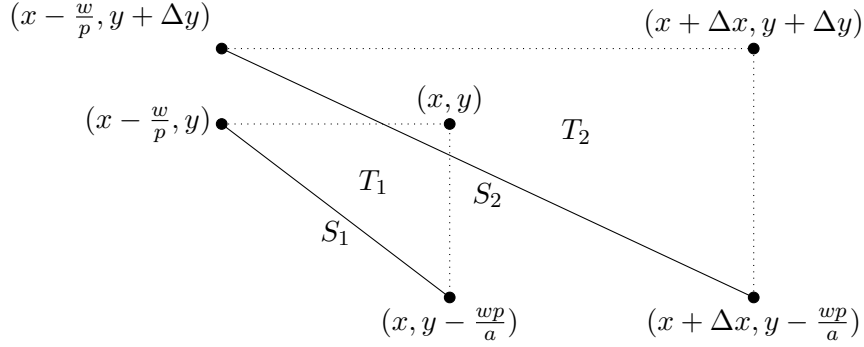
We claim that each point $(x',y')\in A$
is to the south-west of $S_2$ (i.e.~$(x'',y'')\in S_2$ for some
$x''\ge x'$, $y''\ge y'$). If $0\le x'\le x$, $0\le y'\le y$ and
$u(x',y')\le z$, we have
$\smoothz{u}{z}(x,y)=\smoothzp{u}{z}{p}(x,y)\le f(x',y')$ by the
definition of $\smoothzp{u}{z}{p}(x,y)$. Since $f$ is strictly
decreasing in
both variables and takes the value $\smoothz{u}{z}(x,y)$ on $S_1$,
the point $(x',y')$ must therefore be to the south-west of $S_1$.
Thus, $u(x',y')>z$ for every point $(x',y')$ in the (possibly empty)
triangle
\[
T_1:=\{(x',y')\ :\ x'\le x,\ y'\le y,\ (x',y')\text{\ strictly north-east of\ }S_1\};
\]
see \cref{fig:triangles}.
Since each point $(x',y')$ in the triangle
\[
T_2:=\{(x',y')\ :\ x'\le x+\Deltax,\ y'\le y+\Deltay,\ (x',y')\text{\ strictly north-east of\ }S_2\}
\]
is north-east of some point in $T_1\cup\{(x,y)\}$, the assumption that
$u(x,y)>z$ implies that
$u(x',y')>z$ for each $(x',y')$ in $T_2$ too. We conclude that
$A$ and $T_2$ are disjoint, and hence each point in $A$ is to the
south-west of $S_2$ as we claimed.

By the inequality of arithmetic and geometric means, we have
\begin{multline*}
\bigl[g(x+\Deltax,y-\tfrac{wp}{a})-z\bigr]^2=\bigl[\tfrac{a}{q}(\Deltay+\tfrac{wp}{a})\bigr]^2
=\{\text{def.~of }q\}=(w+p\Deltax)(w+\tfrac{a}{p}\Deltay)\\
=w^2+a\Deltax\Deltay+(p\Deltax+\tfrac{a}{p}\Deltay)w
\ge w^2+a\Deltax\Deltay+2w\sqrt{a\Deltax\Deltay}
\\=\left(w+\sqrt{a\Deltax\Deltay}\right)^2,
\end{multline*}
and analogously,
\begin{multline*}
\bigl[g(x-\tfrac{w}{p},y+\Deltay)-z\bigr]^2=\bigl[q(\Deltax+\tfrac{w}{p})\bigr]^2
=(w+p\Deltax)(w+\tfrac{a}{p}\Deltay)\\
\ge\left(w+\sqrt{a\Deltax\Deltay}\right)^2.
\end{multline*}
Since $g$ is an affine function, it follows that
$g(x',y')\ge z+w+\sqrt{a\Deltax\Deltay}
=\smoothz{u}{z}(x,y)+\sqrt{a\Deltax\Deltay}$ for each $(x',y')\in S_2$
and thus for each $(x',y')$ south-west of $S_2$ since $g$ is decreasing
in both variables. Since all of $A$ is south-west of $S_2$ we are done.
\end{proof}

In \cref{lm:largeproduct} below we will use \cref{lm:largedifference} to bound
$\ddxddy{\smooth{u}}$ from below, but first we need a small technical lemma.
\begin{lemma}\label{lm:smoothz}
For any $Z>\smooth{u}(x,y)$ it holds that
\[
\smooth{u}(x,y)=\inf_{0\le z<Z}\smoothz{u}{z}(x,y).
\]
\end{lemma}
\begin{proof}
By definition,
$\smooth{u}(x,y)=\inf_{z\ge 0}\smoothz{u}{z}(x,y)$,
which equals
$\inf_{0\le z<Z}\smoothz{u}{z}(x,y)$ or $\inf_{z\ge Z}\smoothz{u}{z}(x,y)$.
But it is not equal to the last term,
since, by \cref{lm:basicsmoothingA}, $\smoothz{u}{z}(x,y)\ge z$, so
$\inf_{z\ge Z}\smoothz{u}{z}(x,y) \ge Z > \smooth{u}(x,y)$.
\end{proof}

\begin{lemma}\label{lm:largeproduct}
Let $(x,y)$ be a point where $u$ and $\smooth{u}$ are both differentiable
and $\smooth{u}(x,y)<u(x,y)$. Then $\ddxddy{\smooth{u}}\ge a/4$ at $(x,y)$.
\end{lemma}
\begin{proof}
Let $\alpha:=\frac{\partial\smooth{u}}{\partial x}(x,y)$ and
$\beta:=\frac{\partial\smooth{u}}{\partial y}(x,y)$.
Since $\smooth{u}$ is continuous at $(x,y)$, there is a $\delta>0$ such that
$\smooth{u}(x+\delta,y+\delta)<u(x,y)$, and from \cref{lm:smoothz} it follows that
\begin{align}
\smooth{u}(x,y)&=\inf_{0\le z<u(x,y)}\smoothz{u}{z}(x,y)
\text{\ and}\label{eq:firstinf}\\
\smooth{u}(x+\Deltax,y+\Deltay)&=\inf_{0\le z<u(x,y)}\smoothz{u}{z}(x+\Deltax,y+\Deltay)
\label{eq:secondinf}
\end{align}
for any positive $\Deltax$ and $\Deltay$ smaller than $\delta$.
By \cref{lm:largedifference}, for any $z<u(x,y)$ and any
$\Deltax,\Deltay>0$, it holds that
\[
\smoothz{u}{z}(x+\Deltax,y+\Deltay)-\smoothz{u}{z}(x,y)\ge\sqrt{a\Deltax\Deltay}.
\]
Combining this with \cref{eq:firstinf,eq:secondinf} and then applying
the first inequality of \cref{lm:infdiffbounds} yields that
\[
\smooth{u}(x+\Deltax,y+\Deltay)-\smooth{u}(x,y)\ge\sqrt{a\Deltax\Deltay}
\]
for any $0<\Deltax,\Deltay<\delta$.

For any $t>0$,
the directional derivative of
$\smooth{u}$ along the vector
$\bigl(\begin{smallmatrix}t\\ t^{-1}\end{smallmatrix}\bigr)$
at the point $(x,y)$ is $\alpha t +\beta t^{-1}$, so
\[
\alpha t +\beta t^{-1}=\lim_{\eps\rightarrow 0^+}
\frac{\smooth{u}(x+\eps t,y+\eps t^{-1})-\smooth{u}(x,y)}{\eps}\ge\sqrt{a}.
\]
Letting $t$ tend to infinity shows that $\alpha>0$ and letting $t$ tend to
zero shows that $\beta>0$. Finally, letting $t=\sqrt{\beta/\alpha}$ yields
$\alpha\beta\ge a/4$.
\end{proof}

Finally, we have all the tools we need to prove \cref{th:largeproduct}.
\largeproducttheorem*
\begin{proof}
Choose $a:=4C$ and construct $\smooth{u}$ as in
\cref{df:smoothu} with respect to this $a$.
By \cref{lm:basicsmoothingC}, $\smooth{u}$ is doubly increasing,
and by \cref{lm:basicsmoothingB},
$0\le\smooth{u}\le u$ everywhere.
The continuity of $\smooth{u}$ follows from \cref{lm:smoothcontinuous},
and the final property follows from \cref{lm:largeproduct}.
\end{proof}

\section{\texorpdfstring{$\functional_\rho$}{Frho} has a continuous maximizer}\label{sec:mainproof}
In this section we prove \cref{th:main}.
The idea is to take any doubly increasing function $u$ and use \cref{th:largeproduct} to construct
a continuous function $\smooth{u}$ such that $\functional_\rho(\smooth{u})\ge \functional_\rho(u)$.

First, we note that a doubly increasing function can be extended to
any larger domain without increasing the diameter of its image. This was shown
in \cite[Lemma~9.1]{SjostrandIncreasing23}, but the proof is so short that we
repeat it here:
\begin{lemma}\label{lm:increasingextension}
Let $A\subseteq B\subseteq\bbR^2$.
For any $u\in\calU(A)$ there is a $w\in\calU(B)$ such that
\begin{itemize}
\item
the restriction of $w$ to $A$ is $u$,
\item the images $u(A)$ and $w(B)$ have the same closure, and
\item $w(x,y)=\inf_A u$ at any point $(x,y)\in B$ such that $x'>x$ or $y'>y$ for any point $(x',y')\in A$.
\end{itemize}
\end{lemma}
\begin{proof}
For each $(x,y)\in B$, let $P(x,y):=\{(x',y')\in A\,:\,x'\le x,\,y'\le y\}$.
Define $w$ by letting $w(x,y):=\sup_{P(x,y)} u$ if $P(x,y)$ is nonempty, and
$w(x,y):=\inf_A u$ if $P(x,y)$ is empty.
\end{proof}
This enables us to prove the following corollary of \cref{th:largeproduct},
which will be suitable to our application.
\begin{coro}\label{cor:largeproduct}
Let $\Omega$ be a bounded open subset of $\bbR^2$ and let
$u$ be a bounded doubly increasing function on $\Omega$. Then, for any $C>0$
there is
a continuous doubly increasing function $\smooth{u}$ on $\Omega$
such that $\diam\smooth{u}(\Omega)\le\diam u(\Omega)$ and,
for any point $(x,y)\in\Omega$ where $u$ and $\smooth{u}$
are both differentiable, the following holds.
\begin{thmlist}
\item\label{it:coincide}
If $\smooth{u}$ and $u$ coincide at $(x,y)$, their partial derivatives coincide there too.
\item\label{it:otherwise}
Otherwise, $\ddxddy{\smooth{u}}\ge C$ at $(x,y)$.
\end{thmlist}
\end{coro}
\begin{proof}
Without loss of generality, we may assume that $\Omega$ is contained
in the open first quadrant and that $\inf_\Omega u = 0$.
By \cref{lm:increasingextension} we can extend $u$ to the closed first quadrant
$\bbR_{\ge0}^2$ without extending the closure of its image, and so that $u(t,0)=u(0,t)=0$ for any $t\ge0$. Now the corollary follows immediately from \cref{th:largeproduct}, except for property \reflocal{it:coincide}.
But this property is a consequence of the inequality $\smooth{u}\le u$ holding
everywhere, since the nonnegative function $u-\smooth{u}$ has a local minimum
at any point where $u$ and $\smooth{u}$ coincide.
\end{proof}

\begin{prop}\label{pr:utilde}
Let $(\Omega,\rho)$ be a bounded probability domain.
Then, for any $u\in\calU(\Omega)$ there is a
continuous $\smooth{u}\in\calU(\Omega)$ such that
$\diam\smooth{u}(\Omega)\le\diam u(\Omega)$ and
$\functional_\rho(\smooth{u})\ge\functional_\rho(u)$.
\end{prop}
\begin{proof}
Choose $C:=\sup_\Omega\rho$ and choose a $\smooth{u}$ according to \cref{cor:largeproduct}.

Consider a point $(x,y)$ where $u$ and $\smooth{u}$ are both differentiable.
We claim that
\begin{equation}\label{eq:Lclaim}
\textstyle
L(\rho,\ddxddy{\smooth{u}})\ge L(\rho,\ddxddy{u})
\end{equation}
at $(x,y)$. From this claim it follows that $\functional_\rho(\smooth{u})\ge\functional_\rho(u)$.

If $\rho(x,y)=0$ our claim holds, so we may assume that $\rho(x,y)>0$.
By \cref{cor:largeproduct}, either the partial derivatives of $\smooth{u}$ and $u$
coincide at $(x,y)$, in which case our claim holds trivially, or
$\ddxddy{\smooth{u}}\ge C$ at $(x,y)$. In the latter case,
$\sqrt{2\ddxddy{\smooth{u}}/\rho}\ge\sqrt2$ and hence, by \cref{pr:phioneisone},
$\Phi(\sqrt{2\ddxddy{\smooth{u}}/\rho})=1$ and
$L(\rho,\ddxddy{\smooth{u}})=\rho$
at $(x,y)$, so our claim holds also in this case.
\end{proof}

Now, our main theorem follows easily.
\maintheorem*
\begin{proof}
By \cref{th:limitsurface} there exists a maximizer $u$ of $\functional_\rho$ over
$\calU_r(\Omega)$, and by \cref{pr:utilde} there is a continuous
$\smooth{u}\in\calU_r(\Omega)$ such that $\functional_\rho(\smooth{u})\ge \functional_\rho(u)$.
Since $u$ is a maximizer, we have $\functional_\rho(\smooth{u})=\functional_\rho(u)$ and
it follows that $\smooth{u}$ is a maximizer as well.
\end{proof}

\section{Can we find a differentiable maximizer of \texorpdfstring{$\functional_\rho$}{Frho}?}\label{sec:nondifferentiability}
One might ask if it is possible to replace continuity in \cref{th:main} by differentiability. In this section we show that all maximizers might be nondifferentiable, and in that
sense \cref{th:main} is the best possible result.

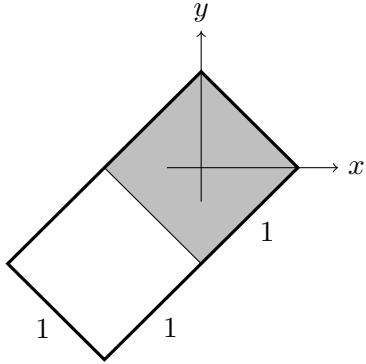
\begin{figure}
\begin{center}
\begin{tikzpicture}[scale=0.9]
\draw[thin,fill=lightgray,rotate around={-45:(0,0)}] (-1,-1) rectangle (1,1);
\draw[very thick,rotate around={-45:(0,0)}] (-1,-3) rectangle (1,1);
\draw [rotate around={-45:(0,0)}]
(1,1) -- (1,-1) node [black,midway,right,xshift=0pt,yshift=-6pt]
{1};
\draw [rotate around={-45:(0,0)}]
(1,-1) -- (1,-3) node [black,midway,right,xshift=0pt,yshift=-6pt]
{1};
\draw [rotate around={-45:(0,0)}]
(1,-3) -- (-1,-3) node [black,midway,below,xshift=-5pt,yshift=1pt]
{1};
\draw [->,thin] (0,-0.5) |- (0,2) node [above] {$y$};
\draw [->,thin] (-0.5,0) |- (2,0) node [right] {$x$};
\end{tikzpicture}
\end{center}
\caption{The rectangle $\Omega$ (thick border) containing the diamond $\Omega'$ (shaded).}
\label{fig:nondifferentiabilityexample}
\end{figure}

Let $\Omega$ be the open rectangle
\[
-\sqrt2<x+y<1/\sqrt2,\ \abs{x-y}<1/\sqrt2,
\]
let $\Omega'\subset\Omega$ be the open diamond $\abs{x}+\abs{y}<1/\sqrt2$,
and let
\[
\rho(x,y):=\begin{cases}
1 & \text{if $(x,y)\in\Omega'$,}\\
0 & \text{if $(x,y)\not\in\Omega'$.}
\end{cases}
\]
See \cref{fig:nondifferentiabilityexample} for an illustration.

Suppose $u$ is a maximizer of $\functional_\rho$ over $\calU_{\sqrt2}(\Omega)$. Then the restriction of $u$ to $\Omega'$ must be a maximizer
of $\functional_1$ over $\calU_{\sqrt2}(\Omega')$. It follows from the proof\footnote{The proposition assumes that $\Phi$ is strictly concave on $[0,\sqrt2]$, but we need only the special case where $c=\sqrt{\rho/ab}$ and for this case the proof works even without the assumption.} of~\cite[Prop.~12.7]{SjostrandIncreasing23} that the function $u_{\rm linear}(x,y):=x+y$ is the unique maximizer of $\functional_1$ over $\calU_{\sqrt2}(\Omega')$ up to an additive constant, so we conclude that
\[
u(x,y)=\begin{cases}
x+y & \text{if $(x,y)\in\Omega'$,}\\
-1/\sqrt{2} & \text{if $(x,y)\not\in\Omega'$.}
\end{cases}
\]
Clearly, $u$ is not differentiable at the boundary between $\Omega'$ and
$\Omega\setminus\Omega'$, that is, at points $(x,y)$ where $x+y=-1/\sqrt2$.

\section{Open problems}\label{sec:future}
As we saw in \cref{sec:nondifferentiability}, continuity cannot be replaced by
differentiability in \cref{th:main}, but the counterexample had a discontinuous
$\rho$.
We leave as an open question whether continuity
or even higher regularity of $\rho$ would guarantee the existence of
a differentiable or even more regular maximizer of $\functional_\rho$.

Another natural question is whether \cref{th:largeproduct} generalizes to higher
dimensions, that is, to increasing functions of more than two variables.
It is known (see \cite{ChabrillacCrouzeix87}) that such functions are differentiable
almost everywhere, and it seems feasible that the ideas in our proof of \cref{th:largeproduct} could be useful also in higher dimension.
\begin{quest}
Let $u\,:\,\bbR_{\ge0}^m\rightarrow\bbR_{\ge0}$ be a function increasing in each
variable and suppose $u(x_1,\dotsc,x_m)=0$ if at least one coordinate is zero.

For any $C>0$, is there
a continuous function $\smooth{u}$ on $\bbR_{\ge0}^m$, increasing in each variable,
such that $0\le \smooth{u}\le u$ everywhere and
$\dd{\smooth{u}}{x_1}\dotsm\dd{\smooth{u}}{x_m}\ge C$ on the set where $\smooth{u}\ne u$
and $\smooth{u}$ is differentiable?
\end{quest}

\section{Acknowledgement}
The author is grateful to an anonymous referee for valuable suggestions.
This work was supported by the Swedish Research Council (reg.no.~2020-04157).

\bibliographystyle{abbrv}
\bibliography{increasing}

\begin{thebibliography}{10}

\bibitem{BassinoEtAl2022}
F.~Bassino, M.~Bouvel, V.~F\'eray, L.~Gerin, M.~Maazoun, and A.~Pierrot.
\newblock Scaling limits of permutation classes with a finite specification: A
  dichotomy.
\newblock {\em Adv.~Math.}, 405:108513, 2022.

\bibitem{BasuEtAl2022}
R.~Basu, S.~Ganguly, A.~Hammond, and M.~Hegde.
\newblock Interlacing and scaling exponents for the geodesic watermelon in last
  passage percolation.
\newblock {\em Commun.~Math.~Phys.}, 393:1241--1309, 2022.

\bibitem{BorgaEtAl2023determinantal}
J.~Borga, C.~Boutillier, V.~F\'eray, and P.-L. M\'eliot.
\newblock A determinantal point process approach to scaling and local limits of
  random {Young} tableaux, 2023.
\newblock Preprint at arXiv:2307.11885.

\bibitem{BorgaEtAl2023}
J.~Borga, S.~Das, S.~Mukherjee, and P.~Winkler.
\newblock Large deviation principle for random permutations.
\newblock {\em Int.~Math.~Res.~Not.~IMRN}, page rnad096, 2023.

\bibitem{BurkillHaslamjones}
J.~C. Burkill and U.~S. Haslam-Jones.
\newblock Notes on the differentiability of functions of two variables.
\newblock {\em Journal of the London Mathematical Society}, s1-7(4):297--305,
  1932.

\bibitem{ChabrillacCrouzeix87}
Y.~Chabrillac and J.-P. Crouzeix.
\newblock {\em Continuity and Differentiability Properties of Monotone Real
  Functions of Several Real Variables}, pages 1--16.
\newblock Springer Berlin Heidelberg, Berlin, Heidelberg, 1987.

\bibitem{CupiniEtAl18}
G.~Cupini, F.~Giannetti, R.~Giova, and A.~Passarelli~di Napoli.
\newblock Regularity results for vectorial minimizers of a class of degenerate
  convex integrals.
\newblock {\em J. Differential Equations}, 265(9):4375--4416, 2018.

\bibitem{DauvergneVirag21}
D.~Dauvergne and B.~Vir\'ag.
\newblock The scaling limit of the longest increasing subsequence, 2021.
\newblock Preprint at arXiv:2104.08210.

\bibitem{DeuschelZeitouni95}
J.-D. Deuschel and O.~Zeitouni.
\newblock Liming curves for i.i.d.~records.
\newblock {\em Ann.\ Probab.}, 23(2):852--878, 1995.

\bibitem{Dubach2023second}
V.~Dubach.
\newblock Increasing subsequences of linear size in random permutations and the
  {Robinson}--{Schensted} tableaux of permutons, 2023.
\newblock Preprint at arXiv:2307.05768.

\bibitem{Dubach2023first}
V.~Dubach.
\newblock Locally uniform random permutations with large increasing
  subsequences, 2023.
\newblock Preprint at arXiv:2301.07658.

\bibitem{GiaquintaGiusti82}
M.~Giaquinta and E.~Giusti.
\newblock On the regularity of the minima of variational integrals.
\newblock {\em Acta Math.}, 148:31--46, 1982.

\bibitem{GlebovEtAl2015}
R.~Glebov, A.~Grzesik, T.~Klimo\v{s}ov\'a, and D.~Kr\'al.
\newblock Finitely forcible graphons and permutons.
\newblock {\em J.~Combin.~Theory~Ser.~B}, 110:112--135, 2015.

\bibitem{Greene}
C.~Greene.
\newblock An extension of {Schensted's} theorem.
\newblock {\em Adv.~Math.}, 14:254--265, 1974.

\bibitem{HoppenEtAl}
C.~Hoppen, Y.~Kohayakawa, C.~G. Moreira, B.~R{\'a}th, and R.~M. Sampaio.
\newblock Limits of permutation sequences.
\newblock {\em J.~Combin.~Theory~Ser.~B}, 103:93--113, 2013.

\bibitem{LadyzhenskayaUraltsevaBook68}
O.~A. Ladyhzenskaya and N.~N. Uraltseva.
\newblock {\em Linear and Quasilinear Elliptic Equations}.
\newblock Academic Press, New York and London, 1968.

\bibitem{LoganShepp}
B.~F. Logan and L.~A. Shepp.
\newblock A variational problem for random {Young} tableaux.
\newblock {\em Advances in Math.}, 26(2):206--222, 1977.

\bibitem{MorreyBook68}
C.~B. Morrey.
\newblock {\em Multiple Integrals in the Calculus of Variations}.
\newblock Springer, Berlin, 1968.

\bibitem{Mukherjee2015}
S.~Mukherjee.
\newblock Fixed points and cycle structure of random permutations.
\newblock {\em Electron.~J.~Prob.}, 21:1--18, 2015.

\bibitem{RomikBook}
D.~Romik.
\newblock {\em The Surprising Mathematics of Longest Increasing Subsequences}.
\newblock Cambridge University Press, 2015.

\bibitem{SjostrandIncreasing23}
J.~Sj{\"o}strand.
\newblock Monotone subsequences in locally uniform random permutations.
\newblock {\em Ann.~Probab.}, 41:1502--1547, 2023.

\bibitem{Ulam}
S.~Ulam.
\newblock {Monte Carlo} calculations in problems of mathematical physics.
\newblock In E.~F. Beckenbach, editor, {\em Modern Mathematics for the
  Engineer, Second Series}, pages 261--281. McGraw-Hill, 1961.

\bibitem{VershikKerov}
A.~M. Vershik and S.~V. Kerov.
\newblock Asymptotic behavior of the {Plancherel} measure of the symmetric
  group and the limit form of {Young} tableaux.
\newblock {\em Dokl.~Akad.~Nauk~SSSR}, 233(6):1024--1027, 1977.

\end{thebibliography}

\end{document}